\documentclass[11pt]{article}
\usepackage{amssymb,amsfonts}
\usepackage{amscd}
\usepackage{amsfonts}
\usepackage{amsmath}
\usepackage{amssymb}
\usepackage{amsthm}
\usepackage{bbm}
\usepackage{CJK}
\usepackage{fancyhdr}
\usepackage{graphicx}
\usepackage{hyperref}
\usepackage{indentfirst}
\usepackage{latexsym}
\usepackage{mathrsfs}
\usepackage{xypic}

\usepackage[top=1in,bottom=1in,left=1.25in,right=1.25in]{geometry}
\def\<{\langle}
\def\>{\rangle}
\def\a{\alpha}
\def\b{\beta}
\def\ci{\circ}
\def\c{\cdot}

\def\D{\Delta}

\def\i{\iota}

\def\r{\rho}
\def\lr{\longrightarrow}

\def\o{\otimes}

\def\vp{\varphi}

\def\<{\langle}
\def\>{\rangle}

\date{}
\begin{document}
\renewcommand{\baselinestretch}{1.2}
\renewcommand{\arraystretch}{1.0}
\title{\bf Constructing New Braided $T$-categories over Monoidal
  Hom-Hopf Algebras}
\date{}
\author {{\bf Miman You \quad  Shuanhong Wang \footnote {Corresponding author:  Shuanhong Wang, shuanhwang2002@yahoo.com}}\\
{\small Department of Mathematics, Southeast University}\\
{\small Nanjing, Jiangsu 210096, P. R. of China}}
 \maketitle
\begin{center}
\begin{minipage}{12.cm}

 \noindent{\bf Abstract} Let  ${\sl Aut}_{mHH}(H)$ denote the set
  of all automorphisms of a monoidal Hopf algebra $H$
  with bijective antipode in the sense of Caenepeel and Goyvaerts
  \cite{CG2011} and let $G$ be a crossed
  product group ${\sl Aut}_{mHH}(H)\times {\sl Aut}_{mHH}(H)$.
  The main aim of this paper is to provide new examples
  of braided $T$-category in the sense of Turaev \cite{T2008}.
  For this purpose,  we first introduce a class of new
  categories $ _{H}\mathcal {MHYD}^{H}(A, B)$ of
  $(A, B)$-Yetter-Drinfeld Hom-modules
  with $A , B \in {\sl Aut}_{mHH}(H)$.
   Then we construct a category ${\cal MHYD}(H)
   =\{{}_{H}\mathcal {MHYD}^{H}(A, B)\}_{(A , B )\in G}$
    and show that such category  forms a new braided $T$-category,
    generalizing the main constructions  by Panaite and Staic  \cite{PS2007}.
    Finally we compute an explicit new example of such braided $T$-categories.

\end{minipage}
\end{center}
\section*{0. INTRODUCTION}

Braided $T$-categories introduced by Turaev \cite{T2008} are
of interest due to their applications in homotopy quantum field theories,
 which are generalizations of ordinary topological quantum field theories.
 As such, they are interesting to different research communities
 in mathematical physics (see \cite{FY1989, K2004, T1994, VA2001, VA2005}).
  Although Yetter-Drinfeld modules over Hopf algebras provide examples of such braided
  $T$-categories, these are rather trivial.
  The wish to obtain more interesting homotopy quantum field theories provides a strong
  motivation to find new examples of braided $T$-categories.
\\

The aim of this article is to construct new examples of braided $T$-categories.
 This is achieved by generalizing an existing construction by Panaite and Staic  \cite{PS2007}
 that twists Yetter-Drinfeld modules over a Hopf algebra $H$ by Hopf algebra automorphisms.
   We will generalize this construction to twisted Yetter-Drinfeld modules over so-called
   monoidal Hom-Hopf algebras, which are Hopf algebras in the Hom-category of a monoidal category
    (see \cite{CG2011}).
   We find a suitable generalisation of the notion of a twisted Yetter-Drinfeld
   module for this setting and obtain a category of
   twisted Yetter-Drinfeld modules that is a braided $T$-category
   in the sense of Turaev \cite{T2008}.
\\

The article is organized as follows.
\\

  We will present the background material in Section 1.
  This section contains the relevant definitions on monoidal Hom-Hopf algebras and braided
  $T$-categories necessary for the understanding of the construction.
  In Section 2, we  define the notion of a Yetter-Drinfeld module over
  a monoidal Hom-Hopf algebra that is twisted by two
  monoidal Hom-Hopf algebra automorphisms
  as well as the notion of a monoidal Hom-entwining structure and
  show how such moinoidal Hom-entwining structures are obtained from automorphisms
  of monoidal Hom-Hopf algebras.
\\

  Section 3 first introduces the tensor product of twisted Yetter-Drinfeld Hom-modules
  and then shows that the twisted Yetter-Drinfeld Hom-modules form a
   braided $T$-category in the sense of Turaev [13]. At the end of the section,
    we give an example of a monoidal Hom-Hopf algebra,
     which can be viewed as a generalization of
     Sweedler's Hopf algebra. And furthermore, we compute an example of a twisted
     Yetter-Drinfeld module over a monoidal Hom-Hopf algebra.
\\

\section*{1. PRELIMINARIES}
\def\theequation{1. \arabic{equation}}
\setcounter{equation} {0} \hskip\parindent

Throughout, let $k$ be a fixed field. Everything is over $k$ unless
 otherwise specified.  We
 refer the readers to the books of Sweedler \cite{S1969}
  for the relevant concepts on the general theory of Hopf
 algebras.  Let $(C, \Delta )$ be a coalgebra. We use the "sigma" notation for
 $\Delta $ as follows:
 $$
 \Delta (c)=\sum c_1\otimes c_2, \,\,\forall c\in C.
 $$

\vskip 0.5cm
 {\bf 1.1. Braided $T$-categories.}
\vskip 0.5cm

 A {\sl monoidal category} ${\cal C}=({\cal C},\mathbb{I},\otimes,a,l,r)$
 is a category ${\cal C}$ endowed with a functor
 $\otimes: {\cal C}\times{\cal C}\rightarrow{\cal C}$
 (the {\sl tensor product}), an object $\mathbb{I}\in {\cal C}$
 (the {\sl tensor unit}), and natural isomorphisms $a$
 (the {\sl associativity constraint}), where
 $a_{U,V,W}:(U\otimes V)\otimes W\rightarrow U\otimes (V\otimes W)$
 for all $U,V,W\in {\cal C}$, and $l$ (the {\sl left unit constraint})
 where $l_U: \mathbb{I}\otimes U\rightarrow U,\,r$
 (the {\sl right unit constraint}) where
 $r_{U}:U\otimes{\cal C}\rightarrow U$ for all $U\in {\cal C}$,
 such that for all $U,V,W,X\in {\cal C},$
 the {\sl associativity pentagon}
 $a_{U,V,W\otimes X}\circ a_{U\otimes V,W,X}
 =(U\otimes a_{V,W,X})\circ a_{U,V\otimes W,X}\circ
 (a_{U,V,W}\otimes X)$ and
 $(U\otimes l_V)\circ(r_U\otimes V)=a_{U,I,V}$ are satisfied.
 A monoidal categoey ${\cal C}$ is {\sl strict} when all
  the constraints are identities.
\\

  Let $G$ be a group and let $Aut({\cal C})$ be the group of
  invertible strict tensor functors from ${\cal C}$ to itself.
   A category ${\cal C}$ over $G$
  is called a {\sl crossed category } if it satisfies the following:
  \begin{eqnarray*}
 &\blacklozenge & {\cal C} \mbox{ is a monoidal category;}\\
 &\blacklozenge & {\cal C} \mbox{ is disjoint union of a family
 of subcategories }\{{\cal C}_{\a }\}_{\a \in
 G},\mbox{ and for any }U\in {\cal C}_{\a },\\
 &&V\in {\cal C}_{\b }, U\o V\in {\cal C}_{\a \b }.
 \mbox{ The subcategory }{\cal C}_{\a }
 \mbox{ is called the }\a\mbox{th component of }{\cal C};\\
 &\blacklozenge & \mbox{Consider a group homomorphism }
   \vp : G\lr Aut({\cal C}), \b \mapsto \vp _{\b }, \mbox{ and
  assume that}\\
  &&  \vp _{\b }(\vp _{\a })
  =\vp _{\b\alpha\beta^{-1}},
  \mbox{ for all }\alpha,\beta\in G.\mbox{ The functors }
  \vp _{\b } \mbox{ are called conjugation}\\
  &&\mbox{ isomorphisms.}
 \end{eqnarray*}

 Furthermore, $ {\cal C}$ is called strict when it is
 strict as a monoidal category.
\\

 {\sl Left index notation}: Given $\a \in G$
 and an object $V\in {\cal C}_{\a }$, the functor $\vp _{\a }$
 will be denoted by ${}^V( \cdot )$, as in Turaev \cite{T2008} or
 Zunino \cite{Z2004}, or even ${}^{\a }( \cdot )$.
 We use the notation ${}^{\overline{V}}( \cdot )$
 for ${}^{\a ^{-1}}( \cdot )$. Then we have
 ${}^V id_U=id_{{} V^U}$ and
 ${}^V(g\circ f)={}^Vg\circ {}^Vf$.
 Since the conjugation $\vp : G\lr Aut({\cal C})$ is a
 group homomorphism, for all $V, W\in {\cal C}$, we have ${}^{V\o W}( \cdot )
 ={}^V({}^W( \cdot ))$ and ${}^\mathbb{I}( \cdot )={}^V({}^{\overline{V}}( \cdot ))
 ={}^{\overline{V}}({}^V( \cdot ))=id_{\cal C}$. Since, for all
 $V\in {\cal C}$, the functor ${}^V( \cdot )$ is strict, we have
 ${}^V(f\o g)={}^Vf\o {}^Vg$, for any morphisms $f$ and $g$ in ${\cal C}$,
  and ${}^V\mathbb{I}=\mathbb{I}$.
\\

 A {\sl braiding} of a crossed category ${\cal C}$ is
 a family of isomorphisms $({c=c_{U,V}})_{U,V}\in {\cal C}$,
 where $c_{U,V}: U\otimes V\rightarrow {}^UV\otimes U$
 satisfying the following conditions:\\
 a) For any arrow $f\in {\cal C}_{\a }(U, U')$ and
 $g\in {\cal C}(V, V')$,
 $$
 (({}^{\a }g)\o f)\circ c _{U, V}=c _{U' V'}\circ (f\o g).
 $$
 b)
 For all $ U, V, W\in {\cal C},$ we have
 $$
 c _{U\o V, W}=a_{{}^{U\o V}W, U, V}\circ (c _{U, {}^VW}\o
 id_V)\circ a^{-1}_{U, {}^VW, V}\circ (\i _U\o c _{V, W})
  \circ a_{U, V, W},
  $$
    $$c _{U, V\o W}=a^{-1}_{{}^UV, {}^UW, U}
 \circ (\i _{({}^UV)}\o c _{U, W})\circ a_{{}^UV, U, W}\ci
 (c _{U, V}\o \i_W)\circ a^{-1}_{U, V, W},$$
 where $a$ is the natural isomorphisms in the tensor category
 ${\cal C}$.\\
 c) For all $ U, V\in {\cal C}$ and $\b\in G$,
 $$ \vp _{\b }(c
 _{U, V})=c _{\vp _{\b }(U), \vp _{\b }(V)}.
 $$

 A crossed category endowed with a braiding is called
 a {\sl braided $T$-category}.
 \\

 {\bf 1.2. Monoidal Hom-Hopf algebras.}
\vskip 0.5cm

Let $\mathcal{M}_{k}=(\mathcal{M}_{k},\o,k,a,l,r )$
 denote the usual monoidal category of $k$-vector spaces and linear maps between them.
 Recall from \cite{CG2011}
 that there is the {\it monoidal Hom-category} $\widetilde{\mathcal{H}}(\mathcal{M}_{k})=
 (\mathcal{H}(\mathcal{M}_{k}),\,\o,\,(k,\,id),
 \,\widetilde{a},\,\widetilde{l},\,\widetilde{r })$, a new monoidal category,
  associated with $\mathcal {M}_{k}$ as follows:

 $\bullet$  The objects of
 $ \mathcal{H}(\mathcal{M}_{k})$ are couples
 $(M,\mu)$, where $M \in \mathcal {M}_{k}$ and $\mu \in Aut_k(M)$, the set of
   all $k$-linear automomorphisms of $M$;

  $\bullet$  The morphism $f:(M,\mu)\rightarrow (N,\nu)$ in $ \mathcal{H}(\mathcal{M}_{k})$
  is the $k$-linear map $f: M\rightarrow N$ in $\mathcal{M}_{k}$
  satisfying   $ \nu \circ f = f\ci \mu$, for any two objects
  $(M,\mu),(N,\nu)\in \mathcal{H}(\mathcal{M}_{k})$;

 $\bullet$   The tensor product is given by
 $$
 (M,\mu)\o (N,\nu)=(M\o N,\mu\o\nu )
 $$
for any $(M,\mu),(N,\nu)\in \mathcal{H}(\mathcal{M}_{k})$.

$\bullet$ The tensor unit is given by $(k, id)$;

$\bullet$   The associativity constraint $\widetilde{a}$
 is given by the formula
 $$
 \widetilde{a}_{M,N,L}=a_{M,N,L}\circ((\mu\o id)\o
  \varsigma^{-1})=(\mu\o(id\o\varsigma^{-1}))\circ a_{M,N,L},
  $$
 for any objects $(M,\mu),(N,\nu),(L,\varsigma)\in \mathcal{H}(\mathcal{M}_{k})$;

 $\bullet$  The left and right unit constraint
  $\widetilde{l}$ and $\widetilde{r }$ are given by
 $$
 \widetilde{l}_M=\mu\circ l_M=l_M\circ(id\o\mu),\, \quad
  \widetilde{r}_M =\mu\circ r_M=r_M\circ(\mu\o id)
  $$
for all $(M,\mu) \in \mathcal{H}(\mathcal{M}_{k})$.
\\

We  now recall from \cite{CG2011} the following notions used later.
\\

 A {\it unital monoidal Hom-associative algebra} (a monoidal Hom-algebra in
 Proposition 2.1 of \cite{CG2011}) is a vector space $A$
 together with an element $1_A\in A$ and linear maps
$$m:A\o A\rightarrow A;\,\,a\o b\mapsto ab, \,\,\,\alpha\in Aut_k(A)$$
such that
\begin{equation}
\alpha(a)(bc)=(ab)\alpha(c),
  \end{equation}
$$\alpha(ab)=\alpha(a)\alpha(b),$$
  \begin{equation}
  a1_A=1_Aa=\alpha(a),
  \end{equation}
\begin{equation}
\alpha(1_A)=1_A ,
  \end{equation}
for all $a,b,c\in A.$
\\

{\bf Remark 1.1.} (1) In the language of Hopf algebras, $m$ is called the
 Hom-multiplication, $\alpha$ is the twisting
 automorphism and $1_A$ is the unit. Note that Eq.(1.1) can be rewirtten as
 $a(b\alpha^{-1}(c)) = (\alpha^{-1}(a)b)c$.
  The monoidal Hom-algebra $A$ with $\alpha$ will be denoted by $(A,\alpha)$.

 (2)  Let $(A,\alpha)$ and  $(A',\alpha')$ be two monoidal Hom-algebras.
   A monoidal Hom-algebra map
 $f:(A,\alpha)\rightarrow (A',\alpha')$ is a linear map such that
 $f\circ \alpha=\alpha'\circ f,f(ab)=f(a)f(b)$ and
 $f(1_A)=1_{A'}.$

(3) The definition of monoidal Hom-algebras is different
 from the unital Hom-associative algebras in \cite{MS2010} and
 \cite{MS2009} in the following sense.
  The same twisted associativity condition (1.1) holds in both cases.
 However, the unitality condition in their notion is the usual untwisted one:
  $a1_A=1_Aa =a,$ for any $a\in A,$ and the twisting map $\alpha$
  does not need to be monoidal (that is, (1.2) and (1.3) are not required).
\\

  {\it A counital monoidal Hom-coassociative coalgebra} is
  an object $(C,\gamma)$ in the category
  $\tilde{\mathcal{H}}(\mathcal{M}_{k})$
 together with linear maps
 $\D:C\rightarrow C\o C,\,\D(c)=c_1\o c_2$ and
 $\varepsilon:C\rightarrow k$ such that
 \begin{equation}
\gamma^{-1}(c_1)\o\D(c_2)=\D(c_1)\o\gamma^{-1}(c_2),
  \end{equation}
 \begin{equation}
\D(\gamma(c))=\gamma(c_1)\o\gamma(c_2),
  \end{equation}
 $$c_1\varepsilon(c_2)=\gamma^{-1}(c)=\varepsilon(c_1)c_2,$$
 \begin{equation}
 \varepsilon(\gamma(c))=\varepsilon(c)
  \end{equation}
 for all $c\in C.$
\\

 {\bf Remark 1.2.} (1) Note that (1.4)is equivalent to
 $c_1\o c_{21}\o \gamma(c_{22})=\gamma(c_{11})\o c_{12}\o c_2.$
 Analogue to monoidal Hom-algebras, monoidal Hom-coalgebras
 will be short for counital monoidal Hom-coassociative coalgebras
 without any confusion.

 (2)  Let $(C,\gamma)$ and $(C',\gamma')$ be two monoidal Hom-coalgebras.
 A monoidal Hom-coalgebra map $f:(C,\gamma)\rightarrow(C',\gamma')$
 is a linear map such that $f\circ \gamma=\gamma'\circ f, \D\circ f=(f\o f)\circ\D$
  and $\varepsilon'\circ f=\varepsilon.$
\\

 {\it A monoidal Hom-bialgebra} $H=(H,\alpha,m,1_H,\D,\varepsilon)$
 is a bialgebra in the monoidal category
 $ \tilde{\mathcal{H}}(\mathcal {M}_{k}).$
 This means that $(H,\alpha,m,1_H)$ is a monoidal Hom-algebra and
 $(H,\alpha,\D,\varepsilon)$ is
 a monoidal Hom-coalgebra such that $\D$ and $\varepsilon$
  are morphisms of algebras,
 that is, for all $h,g\in H,$
 $$\D(hg)=\D(h)\D(g),\, \,\, \D(1_H)=1_H\o1_H,\,\,\,\,\,\,
  \varepsilon(hg)=\varepsilon(h)\varepsilon(g), \,\,\,\,\,\varepsilon(1_H)=1.$$
\\

  A monoidal Hom-bialgebra $(H,\alpha)$ is called {\it a monoidal Hom-Hopf algebra}
 if there exists a morphism (called antipode)
 $S: H\rightarrow H$ in $ \tilde{\mathcal{H}}(\mathcal {M}_{k})$
 (i.e., $S\ci \alpha=\alpha\ci S$),
 which is the convolution inverse of the identity morphism $id_H$
 (i.e., $ S*id=1_H\ci \varepsilon=id*S$). Explicitly,  for all $h\in H$,
 $$
 S(h_1)h_2=\varepsilon(h)1_H=h_1S(h_2).
 $$
\\

{\bf Remark 1.3.} (1) Note that a monoidal Hom-Hopf algebra is
 by definition a Hopf algebra in $ \tilde{\mathcal{H}}(\mathcal {M}_{k})$.

 (2)  Furthermore, the antipode of monoidal Hom-Hopf algebras has
   almost all the properties of antipode of Hopf algebras such as
 $$S(hg)=S(g)S(h),\,\,\,\, S(1_H)=1_H,\,\,\,\,
  \D(S(h))=S(h_2)\o S(h_1),\,\,\,\,\,\,\varepsilon\ci S=\varepsilon.$$
 That is, $S$ is a monoidal Hom-anti-(co)algebra homomorphism.
  Since $\alpha$ is bijective and commutes with $S$,
  we can also have that the inverse $\alpha^{-1}$ commutes with $S$,
  that is, $S\ci \alpha^{-1}= \alpha^{-1}\ci S.$
 \\

 In the following, we recall the notions of actions on monoidal
 Hom-algebras and coactions on monoidal Hom-coalgebras.
\\

  Let $(A,\alpha)$ be a monoidal Hom-algebra.
{\it A left $(A,\alpha)$-Hom-module} consists of
 an object $(M,\mu)$ in $\tilde{\mathcal{H}}(\mathcal {M}_{k})$
  together with a morphism
  $\psi:A\o M\rightarrow M,\psi(a\o m)=a\cdot m$ such that
 $$\alpha(a)\c(b\c m)=(ab)\c\mu(m),\,\,
 \,\,\mu(a\c m)=\alpha(a)\c\mu(m),\,\,
 \,\,1_A\c m=\mu(m),$$
 for all $a,b\in A$ and $m \in M.$
\\

 Monoidal Hom-algebra $(A,\alpha)$ can be
 considered as a Hom-module on itself by the Hom-multiplication.
  Let $(M,\mu)$ and $(N,\nu)$ be two left $(A,\alpha)$-Hom-modules.
  A morphism $f:M\rightarrow N$ is called left
   $(A,\alpha)$-linear if
   $f(a\c m)=a\c f(m),f\ci \mu= \nu\ci f.$
  We denoted the category of left $(A,\alpha)$-Hom modules by
  $\tilde{\mathcal{H}}(_{A}\mathcal {M}_{k})$.
\\

 Similarly, let $(C,\gamma)$ be a monoidal Hom-coalgebra.
 {\it A right $(C,\gamma)$-Hom-comodule} is an object
  $(M,\mu)$ in $\tilde{\mathcal{H}}(\mathcal {M}_{k})$
  together with a $k$-linear map
  $\rho_M:M\rightarrow M\o C,\rho_M(m)=m_{(0)}\o m_{(1)}$ such that
 \begin{equation}
 \mu^{-1}(m_{(0)})\o \D_C(m_{(1)})
 =(m_{(0)(0)}\o m_{(0)(1)})\o \gamma^{-1}(m_{(1)}),
  \end{equation}
    \begin{equation}
 \rho_M(\mu(m))=\mu(m_{(0)})\o\gamma(m_{(1)}),
  \ \ \
 m_{(0)}\varepsilon(m_{(1)})=\mu^{-1}(m),
 \end{equation}
for all $m\in M.$
\\

 $(C,\gamma)$ is a Hom-comodule on itself via the Hom-comultiplication.
 Let $(M,\mu)$ and $(N,\nu)$ be two right $(C,\gamma)$-Hom-comodules.
  A morphism $g:M\rightarrow N$ is called right $(C,\gamma)$-colinear
  if $g\ci \mu=\nu\ci g$ and
  $g(m_{(0)})\o m_{(1)}=g(m)_{(0)}\o g(m)_{(1)}.$
  The category of right
  $(C,\gamma)$-Hom-comodules is denoted by
   $\tilde{\cal{H}}(\cal {M}^C)$ .
\\

 Let $(H, \alpha)$ be a monoidal Hom-bialgebra. We now recall from \cite{CZ2014}
 that a monoidal Hom-algebra $(B,\beta)$ is called {\it a left
 $H$-Hom-module algebra}, if $(B,\beta)$ is a left
 $H$-Hom-module with action $\cdot$ obeying the following axioms:
\begin{equation}
h\c(ab)=(h_1\c a)(h_2\c b),\,\,\,\,\,\,\,\,h\c1_B=\varepsilon(h)1_B,
  \end{equation}
for all $a,b\in B,h\in H.$
\\

 Recall from \cite{LS2014} that a monoidal Hom-algebra $(B,\beta)$ is called
 {\it a left $H$-Hom-comodule algebra}, if $(B,\beta)$ is a left
 $H$-Hom-comodule with coaction $\rho$ obeying the following axioms:
 $$\rho(ab) = a_{(-1)}b_{(-1)}\otimes a_{(0)}b_{(0)},
  \ \ \ \rho_l(1_B)= 1_B\otimes 1_H,$$
for all $a,b\in B,h\in H.$
\\

 Let $(H,m,\Delta,\alpha)$ be a monoidal Hom-bialgebra.
 Recall from (\cite{CZ2014,LS2014}) that a {\it left-right Yetter-Drinfeld Hom-module}
 over $(H,\alpha)$ is the object $(M,\cdot,\rho,\mu)$ which
 is both in  $\tilde{\cal{H}}(_{H}\cal {M})$ and $\tilde{\cal{H}}(\cal {M}^{H})$
 obeying the compatibility condition:
 \begin{equation}
 h_{1}\c m_{(0)}\o h_{2}m_{(1)}=(\alpha(h_{2})\c m)_{(0)} \o
 \alpha^{-1}(\alpha(h_{2})\c m)_{(1)})h_{1}.
  \end{equation}

 {\bf Remark 1.4.} (1) The category of all left-right Yetter-Drinfeld Hom-modules
  is denoted by $\tilde{\cal{H}}(_{H}\cal {YD} ^{H})$ with understanding morphism.

 (2) If $(H,\alpha)$ is a monoidal Hom-Hopf algebra with a bijective
 antipode $S$, then the above equality is equivalent to
 $$\r (h\c m)=\alpha(h_{21})\c m_{(0)}\o (h_{22}\a^{-1}(m_{(1)}))S^{-1}(h_{1}), $$
 for all $h\in H$ and $m\in M$.
\\

\section*{2. $(A, B)$-YETTER-DRINFELD HOM-MODULES}
\def\theequation{2. \arabic{equation}}
\setcounter{equation}{0} \hskip\parindent

In this section, we define the notion of a Yetter-Drinfeld module over
  a monoidal Hom-Hopf algebra that is twisted by two monoidal
  Hom-Hopf algebra automorphisms
  as well as the notion of a monoidal Hom-entwining structure and
  show how such moinoidal Hom-entwining structures are obtained
  from automorphisms of monoidal Hom-Hopf algebras.
\\

In what follows, let  $(H,\alpha)$ be a monoidal Hom-Hopf algebra with
the bijective antipode $S$  and let ${\sl Aut}_{mHH}(H)$ denote the set
  of all automorphisms of a monoidal Hopf algebra $H$.
\\

{\bf Definition 2.1.} Let $A, B\in  {\sl Aut}_{mHH}(H)$.
A left-right {\sl $(A, B)$-Yetter-Drinfeld Hom-module} over $(H,\alpha)$
is a vector space $M$ such that:
\begin{eqnarray*}
 &(1)&(M,\cdot,\mu) \mbox{ is a left }H \mbox{-Hom-module;}\\
 &(2)& (M,\rho,\mu) \mbox{ is a right }H \mbox{-Hom-comodule;}\\
 &(3)& \rho \mbox{ and }\cdot \mbox{ satisfy
  the following compatibility condition: }
   \end{eqnarray*}
 \begin{equation}
\r (h\c m)=\alpha(h_{21})\c m_{(0)}\o (B(h_{22})\a^{-1}(m_{(1)}))A(S^{-1}(h_{1})),
\end{equation}
for all $h\in H$ and $m\in M$. We denote by
$ _{H}\mathcal{MHYD}^{H}(A, B)$
the category of left-right $(A, B)$-Yetter-Drinfeld
Hom-modules, morphisms
being $H$-linear $H$-colinear maps.
\\

 {\bf Remark 2.2.} Note that, $A$ and $B$ are bijective,
 Hom-algebra morphisms,
 Hom-coalgebra morphisms,
 and commute with $S$ and $\alpha$.
 \\

{\bf Proposition 2.3.} One has that Eq.(2.1) is equivalent to
 the following equation:
\begin{equation}
h_{1}\c m_{(0)}\o B(h_{2})m_{(1)}=\mu((h_{2}\c\mu^{-1} (m))_{(0)}) \o (h_{2}\c
 \mu^{-1}(m))_{(1)}A(h_{1}).
\end{equation}

{\bf Proof.} Eq.(2.1)$\Longrightarrow$ Eq.(2.2). We compute as follows
\begin{eqnarray*}
 &&\mu((h_{2}\c \mu^{-1}(m))_{(0)}) \o (h_{2}\c\mu^{-1}( m))_{(1)}A(h_{1})\\
&\stackrel{(2.1)}{=}&  \mu(\a(h_{221})\c \mu^{-1}(m)_{(0)})\o(( B(h_{222})\a^{-2}(m_{(1)}))A(S^{-1}(h_{21})))A(h_{1})\\
 &=& \mu(h_{12}\c \mu^{-1}(m)_{(0)})
 \o(( B\a^{-2}(h_{2})\a^{-2}(m_{(1)}))A(S^{-1}\a(h_{112})))A\a^{2}(h_{111})\\
 &=& \mu(\a^{-1}(h_{1}))\c \mu^{-1}(m)_{(0)})\o B(h_{2})m_{(1)}
 =h_{1}\c m_{(0)}\o B(h_{2})m_{(1)}.
\end{eqnarray*}

 For Eq.(2.2) $\Longrightarrow$ Eq.(2.1), we have
\begin{eqnarray*}
&& \alpha(h_{21})\c m_{(0)}\o (B(h_{22})\a^{-1}(m_{(1)}))A(S^{-1}(h_{1}))\\
&\stackrel{(2.2)}{=}& \mu((\alpha(h_{22})\c \mu^{-1}(m))_{(0)})\o \a^{-1}((\a(h_{22})\c\mu^{-1}(m))_{(1)}A\a(h_{21}))AS^{-1}(h_{1})\\
&=& \mu((h_{2})\c \mu^{-1}(m))_{(0)})\o \a^{-1}((h_{2}\c\mu^{-1}(m))_{(1)}A\a(h_{12}))AS^{-1}\a(h_{11})\\
&=& (\a(h_{2})\c m)_{(0)}\o (h_{2}\c \mu^{-1}(m))_{1}(A(h_{12})AS^{-1}(h_{11}))
=(h\c m)_{(0)}\o (h\c m)_{(1)}.
\end{eqnarray*}
This finishes the proof. \hfill $\blacksquare$
\\

{\bf Example 2.4.}
For $A=B=id_{H},$ we have
 $ _{H}\mathcal{MHYD}^{H}(id, id)
 =\mathcal{H}(_{H}\mathcal{YD}^{H})$, the usual monoidal
 Yetter-Drinfeld Hom-module category. (see \cite{CZ2014,LS2014}).
\\

 {\bf Example 2.5.} (1) Take a non-trivial monoidal Hom-Hopf algebra isomorphism
  $B\in {\sl Aut}_{mHH}(H).$ We
 define $(H_B,\alpha)=(H,\alpha)$ as $k$-vector spaces,
and we can consider a right $H$-Hom-comodule structure on
 $H_B$ via the Hom-comultiplication $\Delta$ and
  a left $H$-Hom-module structure on $H_B$ as follows:
$$
h\cdot y=(B(h_2)\alpha^{-1}(y))S^{-1}(\alpha(h_1)).
$$
for all $h\in H,\,y\in H_B.$ Then it is not hard to check that
$H_B\in\,\, _{H}\mathcal{MHYD}^{H}(id, B)$.

 More generally, if $A,\,B\in {\sl Aut}_{mHH}(H),$
 define $H_{(A,B)}$ as follows: $(H_{(A,B)},\alpha)=(H,\alpha)$  as $k$-vector spaces,
 with right $H$-Hom-comodule structure via Hom-comultiplication
and left $H$-Hom-module structure given by:
$$
h\cdot x=(B(h_2)\alpha^{-1}(x))A(S^{-1}(\alpha(h_1))).
$$
for all $h,\,x\in H.$ It is straightforward to check
that $H_{(A,B)}\in\,\, _{H}\mathcal{MHYD}^{H}(A, B)$.
\\

(2) Recall from Example 3.5 in \cite{CWZ2013}
 that $( H_{4} =k\{ 1,\ g,\ x, \ gx\,\},\a , \Delta , \varepsilon , S )$
 is a monoidal Hom-Hopf algebra, where the algebraic structure are given as follows:

$\bullet$   The multiplication
 $"\circ "$ is given by
 \begin{eqnarray*}
 &&1\circ 1=1,\ \ \ \ \ \ \ \ \ \ 1\circ g = g,
 \ \ \ \ \ \ \ \ \ \ \ 1\circ x=cx,
 \ \ \ \ \ \ \ \ \ \  1\circ gx = cgx,\\
 &&g\circ 1 =g,\ \ \ \ \ \ \ \ \ \  g\circ g=1,\ \ \ \ \
 \ \ \ \ \ \ g\circ x=cgx,\ \ \ \ \ \ \ \ \ g\circ gx=cx,\\
 && x\circ 1=cx,
 \ \ \ \ \ \ \ \ x\circ g=-cgx,
 \ \ \ \ \ \ \  x\circ x=0,
 \ \ \ \ \ \ \ \ \ \ x\circ gx=0,\\
 && gx\circ 1=cgx,\ \ \ \ \
  gx\circ g=-cx,\ \ \ \ \ \ \
   gx\circ x=0,\ \ \ \ \ \ \ \ x\circ gx=0;
\end{eqnarray*}

$\bullet$ The automorphism $\a$ is given by
 $$
 \a(1)=1,\ \a(g)=g,\ \a(x)=cx,\ \a(gx)=cgx,
 $$
  for all $0\neq c\in k;$

$\bullet$  The comultiplication $\Delta$ is defined by
 \begin{eqnarray*}
 &&\Delta(1)=1\otimes 1,\ \ \ \ \ \ \ \ \ \
  \ \ \ \ \ \ \ \ \ \ \ \ \ \ \ \ \
  \Delta(g)=g\otimes g,\\
  &&\Delta(x)=c^{-1}(x\otimes 1)+ c^{-1}(g\otimes x),
 \ \ \Delta(gx)=c^{-1}(gx\otimes g)+c^{-1}(1\otimes gx);
 \end{eqnarray*}

$\bullet$ The counit $\varepsilon$ is defined by

 $$
 \varepsilon(1)=1,\ \
 \varepsilon(g)=1,\ \
 \varepsilon(x)=0,\ \  \varepsilon(gx)=0,
 $$
and

$\bullet$ The antipode $S$ is given by

 $$
 S(1)=1,\ \  S(g)=g,\ \   S(x)=-gx,
\ \  S(gx)=-x,
$$
 \\

Still from Example 3.5 in \cite{CWZ2013} that we have the automorphism group of  the monoidal Hom-Hopf algebra
 $H_{4}$: ${\sl Aut}_{mHH}(H_{4})=\{\left(
  \begin{array}{cccc}
    1 & 0 & 0 & 0 \\
    0 & 1 & 0 & 0 \\
    0 & 0 & \lambda & 0 \\
    0 & 0 & 0 & \lambda \\
  \end{array}
\right)\mid 0\neq\lambda\in k\}$
\\
\\

In what follows, we will give an explicit describe on the Yetter-Drinfeld
Hom-modules given in Part (1) above for $H_4$.
\\

   Let $A=\left(
  \begin{array}{cccc}
    1 & 0 & 0 & 0 \\
    0 & 1 & 0 & 0 \\
    0 & 0 & c' & 0 \\
    0 & 0 & 0 & c' \\
  \end{array}
\right)$.
By Part (1), we can define $H_{4 A}=(H_4,\a)$ as $k$-vector spaces,
but with the right $H_4$-Hom-comodule structure
via $\Delta$ and the left $H_4$-module structures as follows:
\begin{eqnarray*}
 &&1\cdot 1=1,\ \ \ \ \ \ \ \ 1\cdot g = g,
 \ \ \ \ \ \ \ 1\cdot x=cx, \ \ \ \ \ \ \ \ 1\cdot gx = cgx,\\
 &&g\cdot 1 =1,\ \ \ \ \ \ \ \ g\cdot g=g,\ \ \ \ \ \ \
 g\cdot x=-cx,\ \ \ \ \ \ g\cdot gx=-cgx,\\
 && x\cdot 1=-c(1+c')gx,
 \ \ \ \ \ \
  x\cdot g=c(1-c')x,\ \ \ \ \ \
   x\cdot x = 0,
 \ \ \ \ \ \
  x\cdot gx=0,\\
 && gx\cdot 1=c(1-c')gx,\ \ \
  \ \ \ gx\cdot g=-c(1+c')x,\ \ \ \
  gx\cdot x=0,
 \ \ \ \ \ x\cdot gx=0,
\end{eqnarray*}
 for any $0\neq c',c\in k$.
 \\

Then one can check that $H_{4 A}\in _{H_4}\mathcal {MHYD}^{H_4}(A,id)$,
 i.e.,  a left-right $(A, id)$-Yetter-Drinfeld Hom-module over $(H_4, \a )$.
 \\

  Let $B=\left(
  \begin{array}{cccc}
    1 & 0 & 0 & 0 \\
    0 & 1 & 0 & 0 \\
    0 & 0 & c'' & 0 \\
    0 & 0 & 0 & c'' \\
  \end{array}
\right)$. Similarly,
  define $(H_{4 B},\a)=(H_4,\a)$  as $k$-vector spaces
  with the right $H_4$-Hom-comodule structure
via $\Delta$ and the left $H_4$-module structures as follows:
 \begin{eqnarray*}
 &&1\cdot 1=1,\ \ \ \ \ \ \ \ \ \ 1\cdot g = g,
 \ \ \ \ \ \ \ \ \ \  1\cdot x=cx,
 \ \ \ \ \ \ \ \ \ \  1\cdot gx = cgx\\
 &&g\cdot 1 =1,\ \ \ \ \ \ \ \ \ \  g\cdot g=g,\ \ \ \ \
 \ \ \ \ \  g\cdot x=-cx,\ \ \ \ \  g\cdot gx=-cgx\\
 && x\cdot 1=-c(1+c'')gx,
 \ \ \ \ \  x\cdot g=c(-1+c'')x,
 \ \ \ \ \ \ \ \ x\cdot x=0,
 \ \ \ \ \ x\cdot gx=0,\\
 && gx\cdot 1=c(-1+c'')gx,\ \ \ \ \
  gx\cdot g=-c(1+c'')x,\ \ \ \ \
   gx\cdot x=0,\ \ \ \ \  x\cdot gx=0,
\end{eqnarray*}
 for any $0\neq c'',c\in k$.

 It is straightforward to see that
 $H_{4 B}\in _{H_4}\mathcal {MHYD}^{H_4}(id,B)$,
  a left-right $(id, B)$-Yetter-Drinfeld Hom-module over $(H_4, \a )$.
 \\

  Let  $A=\left(
  \begin{array}{cccc}
    1 & 0 & 0 & 0 \\
    0 & 1 & 0 & 0 \\
    0 & 0 & c' & 0 \\
    0 & 0 & 0 & c' \\
  \end{array}
\right)$ and $B=\left(
  \begin{array}{cccc}
    1 & 0 & 0 & 0 \\
    0 & 1 & 0 & 0 \\
    0 & 0 & c'' & 0 \\
    0 & 0 & 0 & c'' \\
  \end{array}
\right)$.

  Define $(H_{4 A,B},\a)=(H_4,\a)$  as $k$-vector spaces,
  with the right $H_4$-Hom-comodule structure
 via $\Delta$ and the left $H$-module structures as follows:
 \begin{eqnarray*}
 &&1\cdot 1=1,\ \ \ \ \ \ \ \ \ \ 1\cdot g = g,
 \ \ \ \ \ \ \ \ \ \  1\cdot x=cx,
 \ \ \ \ \ \ \ \ \ \  1\cdot gx = cgx\\
 &&g\cdot 1 =1,\ \ \ \ \ \ \ \ \ \  g\cdot g=g,\ \ \ \ \
 \ \ \ \ \  g\cdot x=-cx,\ \ \ \ \  g\cdot gx=-cgx\\
 && x\cdot 1=-c(c'+c'')gx,
 \ \ \ \ \  x\cdot g=c(-c'+c'')x,
 \ \ \ \ \ \ \ \ x\cdot x=0,
 \ \ \ \ \ x\cdot gx=0,\\
 && gx\cdot 1=c(-c'+c'')gx,\ \ \ \ \
  gx\cdot g=-c(c'+c'')x,\ \ \ \ \
   gx\cdot x=0,\ \ \ \ \  x\cdot gx=0,
 \end{eqnarray*}
 for any $0\neq c'',c\in k$.
 \\

 Then it is straightforward to see that $H_{4 (A,B)}$ is
  a left-right $(A, B)$-Yetter-Drinfeld Hom-module over $(H_4, \a )$,
   i.e.,
 $H_{4 (A,B)}\in _{H_4}\mathcal {MHYD}^{H_4}(A,B).$
 \\

{\bf Definition 2.6.} A left-right {\sl monoidal Hom-entwining structure} is a triple
 $(H, \, C,\, \psi)$, where $(H, \a)$ is a monoidal Hom-algebra
 and $(C,\gamma)$ is a monoidal Hom-coalgebra with a linear map
 $\psi: H\otimes C\rightarrow H\otimes C,
 \,\,\,\,h\otimes c\mapsto _{\psi}h\otimes c^{\psi}$
 satisfying the following conditions:
\begin{equation}
 _{\psi}(hg)\otimes c^{\psi}= _{\phi}h_{\psi}g\otimes \gamma(\gamma^{-1}(c)^{\psi\phi}),
\end{equation}
\begin{equation}
 _{\psi}1\otimes c^{\psi}= 1_A\otimes c,
\end{equation}
\begin{equation}
 _{\psi}h\otimes \Delta(c^{\psi})= \a(_{\phi\psi}\a^{-1}(h))\otimes (c^{\phi}_{1}\otimes c^{\psi}_{2}),
\end{equation}
\begin{equation}
 \varepsilon(c^{\psi}) _{\psi}h= \varepsilon(c)a,
\end{equation}
  Over a monoidal Hom-entwining structure $(H, \, C,\, \psi)$,
 a left-right monoidal entwined Hom-module $M$ is both a right
 $C$-Hom-comodule and a left $H$-Hom-module such that
 $$ \rho^{M}(h\cdot m)=_{\psi}\a^{-1}(h)\cdot m_{(0)}\otimes \a(m_{(1)})^{\psi}$$
 for all $h\in H$ and $m\in M$. We denote the category of all
 monoidal entwined Hom-modules over $(H,\,C,\,\psi)$
 by $_{H}\mathcal{M}^{C}(\psi)$.
 \\

 Let $(H,\a)$ be a monoidal Hom-Hopf algebra with $S$,
 and define a linear map
 $$\psi(A,B):H\otimes H\rightarrow H\otimes H,\
 \ a\otimes c\mapsto _{\psi}a\otimes c^{\psi}
 =\a^{2}(a_{21})\otimes (B(a_{22})\a^{-2}(c))AS^{-1}(a_1),$$
 for all $A, B\in  {\sl Aut}_{mHH}(H)$.
\\

{\bf Proposition 2.7.} With notations as above,
 $(H, H,\psi(A,B))$ is a monoidal Hom-entwining structure
\underline{} for all $A, B\in  {\sl Aut}_{mHH}(H)$.

 {\bf Proof.} We need to prove that Eqs.(2.3-2.6) hold.
 First, it is straightforward to check Eqs.(2.4) and (2.6).
 In what follows, we only verify Eqs.(2.3) and (2.5).
 In fact, for all $a,b,c\in H$, we have
 \begin{eqnarray*}
&&_{\phi}a_{\psi}b\otimes \a(\a^{-1}(c)^{\psi\phi})\\
&=&\a^2(a_{21})_{\psi}b \otimes
\a((B(a_{22})\a^{-2}(\a^{-1}(c)^{\psi}))AS^{-1}(a_1))\\
&=&\a^2(a_{21})\a^2(b_{21})
\otimes\a([B(a_{22})\a^{-2}((B(b_{22})\a^{-3}(c))AS^{-1}(b_1))]AS^{-1}(a_1))\\
&=&\a^2(a_{21}b_{21})\otimes[B\a(a_{22})
((B\a^{-1}(b_{22})\a^{-4}(c))AS^{-1}\a^{-1}(b_1))]AS^{-1}\a(a_1)\\
&=&\a^2(a_{21}b_{21})
\otimes(B(a_{22}b_{22})\a^{-2}(c))(AS^{-1}(b_1)AS^{-1}(a_1))
= _{\psi}(ab)\otimes c^{\psi},
\end{eqnarray*}
and Eq.(2.3) is proven.

For all $a\in H$, we have
\begin{equation}
a_{1}\otimes a_{211}\otimes a_{2121}\otimes a_{2122}\otimes a_{22}
=\a(a_{11})\otimes \a^{-1}(a_{12})\otimes \a^{-2}(a_{21})
\otimes \a^{-1}(a_{221})\otimes \a(a_{222})
\end{equation}

As for Eq.(2.5), we compute:
\begin{eqnarray*}
&& \a(_{\phi\psi}\a^{-1}(a))\otimes (c^{\phi}_{1}\otimes c^{\psi}_{2})\\
&=&\a(\a^{2}((_{\psi}\a^{-1}(a))_{21}))\otimes
((B((_{\psi}\a^{-1}(a))_{22})\a^{-2}(c_{1}))AS^{-1}((_{\psi}\a^{-1}(a))_{1})
\otimes c^{\psi}_{2})\\
&=&\a(\a^{2}(\a^{2}(\a^{-1}(a)_{21})_{21}))\otimes
((B((\a^{2}(\a^{-1}(a)_{21})_{22})\a^{-2}(c_{1}))AS^{-1}((\a^{2}(\a^{-1}(a)_{21})_{1})\\
&&\otimes (B(\a^{-1}(a)_{22})\a^{-2}(c_{2}))AS^{-1}(\a^{-1}(a)_{1}))\\
&=&\a^{4}(a_{2121})\otimes
((B\a(a_{2122})\a^{-2}(c_{1}))AS^{-1}\a(a_{211})
\otimes (B\a^{-1}(a_{22})\a^{-2}(c_{2}))AS^{-1}\a^{-1}(a_{1}))\\
&\stackrel{(2.7)}{=}&\a^{2}(a_{21})\otimes
((B(a_{221})\a^{-2}(c_{1}))AS^{-1}(a_{12})
\otimes (B(a_{222})\a^{-2}(c_{2}))AS^{-1}(a_{11}))\\
&=&\a^{2}(a_{21})\otimes((B(a_{22})\a^{-2}(c))AS^{-1}(a_1))_1
\otimes(B(a_{22})\a^{-2}(c))AS^{-1}(a_1))_2)\\
&=&_{\psi}a\otimes \Delta(c^{\psi}).
\end{eqnarray*}
and Eq.(2.5) is proven.

This finishes the proof. \hfill $\blacksquare$
\\

 By Proposition 2.7, we have a monoidal entwined
  Hom-module category $ _{H}\mathcal{M}^{H}(\psi(A, B))$
 over $(H, H,\psi(A,B))$ with $A,B\in{\sl Aut}_{mHH}(H)$.
 In this case, for all $M\in _{H}\mathcal{M}^{H}(\psi(A, B))$,
 we have
  $$\r (h\c m)=\alpha(h_{21})\c m_{(0)}
  \o (B(h_{22})\a^{-1}(m_{(1)}))A(S^{-1}(h_{1})),$$
  for all $h\in H,m\in M$. Thus means that
  $_{H}\mathcal{M}^{H}(\psi(A, B))=
   _{H}\mathcal{MHYD}^{H}(A, B)$
   as categories.
   \\

  {\bf Definition 2.8.} Let $(H,\alpha)$ be a
  monoidal Hom-algebra. A monoidal Hom-algebra
  $(N,\nu)$ is called an $(H,\alpha)$-Hom-bicomodule
  algebra, if $(N,\nu)$ is a left $(H,\alpha)$-Hom-comodule
  and a right $(H,\alpha)$-Hom-comodule
  with coactions $\rho_r$ and $\rho_l$
  obeying the following axioms:
  \begin{eqnarray*}
 &(1)&\rho_l(n) = n_{[-1]}\otimes n_{[0]},\ \mbox{ and }
  \ \rho_r(n)= n_{<0>}\otimes n_{<1>},\\
 &(2)& (N,\nu)\mbox{ is a left }H \mbox{-Hom-comodule algebra;}\\
 &(3)&  (N,\nu)\mbox{ is a right }H \mbox{-Hom-comodule algebra;}\\
 &(4)& \rho_l \mbox{ and }\rho_l \mbox{ satisfy
  the following compatibility condition: for all }  n\in N,\\
 &\ \ \ \ \ &n_{<0>[-1]}\otimes n_{<0>[0]}\otimes \alpha^{-1}(n_{<1>})
  = \alpha^{-1}(n_{[-1]})\otimes n_{[0]<0>}\otimes n_{[0]<1>}\\
  &\ \ \ \ \ &  \quad \quad =n_{\{-1\}}\otimes n_{\{0\}}\otimes n_{\{1\}}
    \in(H\otimes N)\otimes H = H\otimes( N\otimes H).
 \end{eqnarray*}
\\

   {\bf Example 2.9.} Let $A,B\in {\sl Aut}_{mHH}(H)$,
 and $H_{(A,B)}=H$ as algebra, with $H$-Hom-comodule
  structures as follows:
 for all $h\in H,$
 \begin{eqnarray*}
 && H_{(A,B)}\rightarrow H\otimes H_{(A,B)},\,\,\,
 h\mapsto h_{[-1]}\otimes h_{[0]}=A(h_{1})\otimes h_{2},\\
 && H_{(A,B)}\rightarrow H_{(A,B)}\otimes H,\,\,\,
 h\mapsto h_{<0>}\otimes h_{<1>}=h_{1}\otimes B(h_{2}).
 \end{eqnarray*}
  Then on can check $H_{(A,B)}$ is an $H$-Hom-bimodule algebra.
 \\

 {\bf Definition 2.10.} Let $(H,\alpha)$ be a
 monoidal Hom-Hopf algebra and
 $(N,\nu)$ be an $H$-Hom-bicomodule algebra,
  A left-right Yetter-Drinfeld Hom-module
 is a $k$-modules $(M,\mu)$ together with a left
 $N$-action (denoted by $n\otimes m\mapsto n\cdot m$)
 and a right $H$-coaction (denoted by $ m\mapsto m_{0}\otimes m_{1}$)
 satisfying the eqivalent compatibility conditions:
 \begin{eqnarray*}
 &&(n\cdot m)_{0}\otimes(n\cdot m)_{1}
 =\nu(n_{[0]<0>})\cdot m_{0}
 \otimes( n_{[0]<1>}\alpha^{-1}(m_{1}))S^{-1}(n_{[-1]}),\\
  &&n_{<0>}\cdot m_{0}\otimes n_{<1>} m_{1}
 =\mu((n_{[0]}\cdot \mu^{-1}(m))_{0})
 \otimes(n_{[0]}\cdot \mu^{-1}(m))_{1}n_{[-1]}.
\end{eqnarray*}
 for all $n\in N$ and $m\in M$.
 Then we all $(H,N,H)$ a Yetter-Drinfeld Hom-datum
 $(H,N,H)$(the second $H$ is regarded as an $H$-Hom-bimodule coalgebra).
 Our notion for the category of a left-right Yetter-Drinfeld Hom-modules
 and $N$-linear $H$-colinear maps
 will be  $ _{N}\mathcal{MHYD}^{H}(H).$
  \\

 {\bf Example 2.11.} Let $H_{(A,B)}$ be an $H$-Hom-cobimodule
 algebra, with an $H$-Hom-comodule structures
 shown in Example 2.9.
 Then we can consider the Yetter-Drinfeld Hom-datum
 $(H,H_{(A,B)},H)$ and the Yetter-Drinfeld Hom-modules
 over it, $ _{H_{(A,B)}}\mathcal{MHYD}^{H}(H)$.
 \\

 {\bf Proposition 2.12.}
 $_{H}\mathcal{MHYD}^{H}(A, B)= _{H_{(A,B)}}\mathcal{MHYD}^{H}(H)$.
 \\

 It is easy to see that the compatibility
 conditions for the two categories are the same.
 The easy proof of this is left to the reader.

\section*{3. A BRAIDED $T$-CATEGORY $\mathcal {MHYD}(H)$}
\def\theequation{3. \arabic{equation}}
\setcounter{equation} {0} \hskip\parindent

In this section, we will construct a class of new braided $T$-categories
 $\mathcal {MHYD}(H)$ over any monoidal Hom-Hopf algebra $(H,\alpha)$.
\\

{\bf Proposition 3.1.} If $(M,\mu) \in {_{H}}\mathcal {MHYD}^{H}(A,B)$
 and $(N,\nu)\in {_{H}}\mathcal {MHYD}^{H}(C,D)$, with
 $A,B,C,D\in {\sl Aut}_{mHH}(H)$, then $(M \o N,\mu\o\nu)
 \in {_{H}}\mathcal {MHYD}^{H}(AC, DC^{-1}BC)$ with structures as follows:
\begin{eqnarray*}
&&h\c (m \o n)=C (h_{1})\c m \o C^{-1}BC(h_{2})\c n,\\
&&m\o n \mapsto (m_{(0)}\o n_{(0)})\o n_{(1)}m_{(1)}.
\end{eqnarray*}
for all $m\in M,n\in N$ and $h\in H.$
\\

{\bf Proof.} First, it is easy to get that
$(M \o N,\mu\o\nu)$ is a left $H$-module
and a right $H$-comodule.
Next, we compute the compatibility condition as follows:
\begin{eqnarray*}
&&(h\c (m\o n))_{(0)}\o (h\c (m\o n))_{(1)}\\
 &=&((C(h_{1})\c m)_{(0)} \o (C^{-1}BC(h_{2})\c n)_{(0)})\o
 (C^{-1}BC(h_{2})\c n)_{(1)} C(h_{1}\c m)_{(1)}\\
&\stackrel{(2.1)}{=}&(C\a (h_{121})\c m_{(0)} \o  C^{-1}BC\a(h_{221})\c n_{(0)})\o
[(DC^{-1}BC(h_{222})\a^{-1}(n_{(1)}))\\
&& CS^{-1}C^{-1}BC(h_{21})][(BC(h_{122})\a^{-1}(m_{(1)}))S^{-1}AC(h_{11})]\\
&=&(C(h_{12})\c m_{(0)} \o  C^{-1}BC\a(h_{221})\c n_{(0)})\o
[(DC^{-1}BC(h_{222})\a^{-1}(n_{(1)})) \\
&&S^{-1}BC\a(h_{212})][(BC(h_{211})\a^{-1}(m_{(1)}))S^{-1}AC(h_{11})]\\
&=&(C(h_{12})\c m_{(0)} \o  C^{-1}BC\a(h_{221})\c n_{(0)})\o
(DC^{-1}BC\a(h_{222})n_{(1)}) \\
&&[(BC\a^{-1}(S^{-1}(h_{212})h_{211})\a^{-1}(m_{(1)}))S^{-1}AC(h_{11}))]\\
&=&(C(h_{12})\c m_{(0)} \o  C^{-1}BC\a(h_{221})\c n_{(0)})\o
(DC^{-1}BC\a(h_{222})n_{(1)}) \\
&&[(BC\a^{-1}(\varepsilon(h_{21})1_{H})\a^{-1}(m_{(1)}))S^{-1}AC(h_{11}))]\\
&=&(C(h_{12})\varepsilon(h_{21})\c m_{(0)} \o  C^{-1}BC\a(h_{221})\c n_{(0)})\o
(DC^{-1}BC\a(h_{222})n_{(1)}) \\
&&(m_{(1)}S^{-1}AC(h_{11}))\\
&=&(C\a(h_{121})\varepsilon\a(h_{122})\c m_{(0)} \o  C^{-1}BC(h_{21})\c n_{(0)})\o
(DC^{-1}BC(h_{22})n_{(1)}) \\
&&(m_{(1)}S^{-1}AC(h_{11}))\\
&=&(C(h_{12})\c m_{(0)} \o  C^{-1}BC(h_{21})\c n_{(0)})\o
(DC^{-1}BC(h_{22})(\a^{-1}(n_{(1)}) \a^{-1}(m_{(1)})))\\
&&S^{-1}AC\a(h_{11})\\
&=&(C\a(h_{211})\c m_{(0)} \o  C^{-1}BC\a(h_{212})\c n_{(0)})\o
(DC^{-1}BC(h_{22})(\a^{-1}(n_{(1)}m_{(1)}))\\
&&S^{-1}AC(h_{1})\\
&=&\a(h_{21})\c (m\o n)_{(0)}\o DC^{-1}BC(h_{22})\a^{-1}(m\o
n)_{(1)}AC(S^{-1}(h_{1})).
\end{eqnarray*}
for all $m\in M,n\in N$ and $h\in H.$
This completes the proof.
 \hfill $\blacksquare$
\\

{\bf Remark. 3.2.} Note that, if $(M,\mu)\in {_{H}}\mathcal
{MHYD}^{H}(A,B),\,\, (N,\nu) \in {_{H}}\mathcal {MHYD}^{H}(C,D)$ \\
and
$(P,\varsigma) \in {_{H}}\mathcal {MHYD}^{H}(E, F)$, then\,\, $(M\o
N)\o P=M \o(N\o P)$
as objects \\
in $_{H}\mathcal {MHYD}^{H}(ACE,
FE^{-1}DC^{-1}BCE)$.
\\

 Denote $G={\sl Aut}_{mHH}(H)
 \times {\sl Aut}_{mHH}(H)$
 a group with multiplication as follows:
  for all $A,B,C,D\in {\sl Aut}_{mHH}(H)$,
\begin{equation}
(A,B)\ast (C,D)=(AC, DC^{-1}BC).
\end{equation}
 The unit of this group is $(id,id)$ and $(A,B)^{-1}=(A^{-1},
 AB^{-1}A^{-1})$.

The above proposition means that if
$M \in {_{H}}\mathcal {MHYD}^{H}(A,B)$
 and $N\in {_{H}}\mathcal {MHYD}^{H}(C,D)$, then
 $M \o N \in {_{H}}\mathcal {MHYD}^{H}((A,B)\ast (C,D)).$
\\

 {\bf Proposition 3.3.} Let $(N,\nu) \in {_{H}}\mathcal {MHYD}^{H}(C,D)$
 and $(A,B)\in G$. Define ${}^{(A,B)}N=N$ as
 vector space, with structures: for all $n\in N$ and $h\in H.$
$$h\rhd n=C^{-1}BCA^{-1}(h)\c n,$$
\begin{equation}
n\mapsto n_{<0>}\o n_{<1>}=n_{(0)}\o AB^{-1}(n_{(1)}).
\end{equation}
Then
$${}^{(A,B)}N \in {_{H}}\mathcal {MHYD}^{H}(ACA^{-1},AB^{-1}DC^{-1}BCA^{-1})
={_{H}}\mathcal {MHYD}^{H}((A,B)\ast
(C,D)\ast (A,B)^{-1}).$$

{\bf Proof.} Obviously, the equations above define
a module and a comodule action.
In what follows,
we show the compatibility condition:
\begin{eqnarray*}
&&(h\rhd n)_{<0>}\o(h\rhd n)_{<1>}\\
&=&(C^{-1}BCA^{-1}(h)\c n)_{(0)}\o AB^{-1}((C^{-1}BCA^{-1}(h)\c n)_{(1)})\\
&=&C^{-1}BCA^{-1}\a(h_{21})\c n_{(0)}\o AB^{-1}((DC^{-1}BCA^{-1}(h_{22})\a^{-1}(n_{(1)}))\\
&&CC^{-1}BCA^{-1}S^{-1}(h_{1}))\\
&=&C^{-1}BCA^{-1}\a(h_{21})\c n_{(0)}\o (AB^{-1}DC^{-1}BCA^{-1}(h_{22})AB^{-1}\a^{-1}(n_{(1)}))
ACA^{-1}S^{-1}(h_{1}))\\
&=&\a(h_{21})\rhd n_{<0>}\o (AB^{-1}DC^{-1}BCA^{-1}(h_{22})\a^{-1}(n_{<1>}))
ACA^{-1}S^{-1}(h_{1}))
\end{eqnarray*}
for all $n\in N$ and $h\in H,$
that is ${}^{(A,B)}N \in {_{H}}\mathcal {MHYD}^{H}(ACA^{-1},AB^{-1}DC^{-1}BCA^{-1})$
\hfill $\blacksquare$
\\

{\bf Remark. 3.4.} Let $(M,\mu) \in {_{H}}\mathcal {MHYD}^{H}(A, B),
 \ \ (N,\nu) \in {_{H}}\mathcal
{MHYD}^{H}(C, D),\,\, \mbox {and}\, (E,F)\in G$. Then by the
 above proposition, we have:
 $$
 {}^{(A, B)\ast (E, F)}N={}^{(A, B)}({}^{(E, F)}N),
$$
as objects in $_{H}\mathcal {MHYD}^{H}(AECE^{-1}A^{-1},
AB^{-1}EF^{-1}DC^{-1}FE^{-1}BECE^{-1}A^{-1})$ and
$$
{}^{(E,F)}(M\o N)= {}^{(E,F)}M \o {}^{(E,F)}N,
$$
as objects in $_{H}\mathcal {MHYD}^{H}(EACE^{-1},
 EF^{-1}DC^{-1}BA^{-1}FACE^{-1})$.
\\

 {\bf Proposition 3.5.} Let $(M,\mu) \in {_{H}}\mathcal {MHYD}^{H}(A,B)$
 and $(N,\nu)\in {_{H}}\mathcal {MHYD}^{H}(C,D)$, take
  ${}^{M}N={}^{(A,B)}N$ as explained in Subsection 1.2.
  Define a map $c_{M, N}: M \o N \rightarrow {}^{M}N \o M$ by
\begin{equation}
 c_{M,N}(m\o
 n)=\nu(n_{(0)})\o B^{-1}(n_{(1)})\c \mu^{-1}(m).
\end{equation}
for all $m\in M,n\in N.$
Then $c_{M, N}$ is both an $H$-module map and an $H$-comodule map,
 and satisfies the following formulae
 (for $(P,\varsigma) \in {_{H}}\mathcal {MHYD}^{H}(E,F)$):
\begin{equation}
a^{-1}_{{}^{M\o N}P, M, N}\circ c _{M\o N, P}\circ a^{-1}_{M, N, P}=(c _{M, {}^NP}\o
 id_N)\circ a^{-1}_{M, {}^NP, N}\circ (id _M\o c _{N, P})  ,
\end{equation}
\begin{equation}
 a_{{}^MN, {}^MP, M}
 \circ  c _{M, N\o P}\circ a_{M, N, P}=(id _{{}^MN }\o c _{M, P})\circ a_{{}^MN, M, P}\ci
 (c _{M, N}\o id_P).
\end{equation}

Furthermore, if $(M,\mu) \in {_{H}}\mathcal {MHYD}^{H}(A,B)$ and
$(N,\nu) \in{_{H}}\mathcal {MHYD}^{H}(C,D)$, then\\
 $c_{{}^{(E,F)}M,{}^{(E,F)}N}=c_{M,N},$
 for all $(E,F)\in G$.
\smallskip
\\

{\bf Proof.} First, we prove that $c_{M,N}$ is an $H$-module map.
 Take  $h\cdot(m\o n)=C(h_{1})\c m\o C^{-1}BC(h_{2})\c n$
 and $h\cdot(n\o m)= C^{-1}BC(h_{1})\c n\o B^{-1}DC^{-1}BC(h_{2})\c m$
 as explained in Proposition 3.1.
 \begin{eqnarray*}
&&c_{M,N}(h\cdot(m\o n))\\
&=&\nu((C^{-1}BC(h_{2})\c n)_{(0)})\o
B^{-1}((C^{-1}BC(h_{2})\c n)_{(1)})\c\mu^{-1}(C(h_{1})\c m)\\
&=&\nu(C^{-1}BC\a(h_{221})\c n_{(0)})\o
B^{-1}((DC^{-1}BC(h_{222})\a^{-1}( n_{(1)}))CC^{-1}BCS^{-1}(h_{21}))\\
&&\c\mu^{-1}(C(h_{1})\c m)\\
&=&\nu(C^{-1}BC\a(h_{221})\c n_{(0)})\o
B^{-1}(DC^{-1}BC\a(h_{222}) n_{(1)})\\
&&\c ((CS^{-1}\a^{-1}(h_{21})
C\a^{-2}(h_{1}))\c \mu^{-1}(m))\\
&=&\nu(C^{-1}BC(h_{21})\c n_{(0)})\o
B^{-1}(DC^{-1}BC(h_{22}) n_{(1)})\\
&&\c (C\a^{-1}((S^{-1}(h_{12})
h_{11}))\c \mu^{-1}(m))\\
&=&\nu(C^{-1}BC(h_{21})\c n_{(0)})\o
B^{-1}(DC^{-1}BC(h_{22}) n_{(1)})\c (C\a^{-1}(\varepsilon(h_{1})1_{H})\c \mu^{-1}(m))\\
&=&C^{-1}BC(h_{1})\c \nu (n_{(0)})\o
B^{-1}DC^{-1}BC(h_{2})\c(B^{-1}( n_{(1)})\c \mu^{-1}( m)) \\
&=&\psi_{N\o M}((h\o(\nu(n_{(0)})\o B^{-1}(n_{(1)})\c \mu^{-1}(m)))\\
&=&\psi_{N\o M}\ci(id\o c_{M,N})(h\o(m\o n))
\end{eqnarray*}

Secondly, we check that $c_{M, N}$ is an $H$-comodule map as
 follows:
\begin{eqnarray*}
&&\rho_{N\o M}\ci c_{M,N}(m\o n) \\
&=&((\nu(n_{(0)}))_{<0>}\o( B^{-1}(n_{(1)})\c\mu^{-1}(m))_{(0)})
\o( B^{-1}(n_{(1)})\c\mu^{-1}(m))_{(1)}\\
&&(\nu(n_{(0)}))_{<1>}\\
&=&(\nu(n_{(0)(0)})\o( B^{-1}(n_{(1)})\c\mu^{-1}(m))_{(0)})
\o( B^{-1}(n_{(1)})\c\mu^{-1}(m))_{(1)}\\
&&AB^{-1}\a(n_{(0)(1)})\\
&=&(\nu(n_{(0)(0)})\o B^{-1}\a(n_{(1)21})\c\mu^{-1}(m)_{(0)})
\o((B B^{-1}(n_{(1)22})\a^{-2}(m_{(1)}))\\
&&AB^{-1}S^{-1}(n_{(1)1}))AB^{-1}\a(n_{(0)(1)})\\
&\stackrel{(1.7)}{=}&(n_{(0)}\o B^{-1}\a(n_{(1)21})\c\mu^{-1}(m_{(0)}))
\o(\a(n_{(1)22})\a^{-1}(m_{(1)}))\\
&&(AB^{-1}S^{-1}\a(n_{(1)12})AB^{-1}\a(n_{(1)11}))\\
&=&(n_{(0)}\o B^{-1}\a(n_{(1)21})\c\mu^{-1}(m_{(0)}))
\o(\a^2(n_{(1)22})m_{(1)})\varepsilon(n_{(1)1})\\
&\stackrel{(1.7)}{=}&(\nu( n_{(0)(0)})\o B^{-1}(n_{(0)(1)})\c\mu^{-1}(m_{(0)}))
\o n_{(1)}m_{(1)}\\
&=&(c_{M,N}\o id)((m_{(0)}\o n_{(0)})\o n_{(1)}m_{(1)})
=(c_{M,N}\o id)\rho(m\o n).
\end{eqnarray*}

Then we will check Eqs.(3.4) and (3.5). On the one hand,
\begin{eqnarray*}
&&a^{-1}_{{}^{M\o N}P, M, N}\circ c _{M\o N, P}\circ a^{-1}_{M, N, P}(m\o (n\o p))\\
&=&a^{-1}_{{}^{M\o N}P, M, N}\circ c _{M\o N, P}((\mu^{-1}(m)\o n)\o \varsigma(p))\\
&=&a^{-1}_{{}^{M\o N}P, M, N}(\varsigma^{2} ( p_{(0)} )\o C^{-1}B^{-1}CD^{-1}\a( p_{(1)} )
\c(\mu^{-2}(m)\o \nu^{-1}(n)))\\
&=&(\varsigma ( p_{(0)} )\o  B^{-1}CD^{-1}\a( p_{(1)1})\c \mu^{-2}(m))\o
D^{-1}\a^{2}( p_{(1)2})\c n\\
&=&(\varsigma^{2}(p_{(0)(0)})\o B^{-1}CD^{-1}\a(p_{(0)(1)})\c\mu^{-2}(m))
\o D^{-1}\a(p_{(1)})\c n\\
&=&(\varsigma(\varsigma(p_{(0)})_{<0>})\o B^{-1}(\varsigma(p_{(0)})_{<1>})\c\mu^{-2}(m))
\o D^{-1}\a(p_{(1)})\c n\\
&=&(c _{M, {}^NP}\o id_N)((\mu^{-1}(m)\o \varsigma(p_{(0)}))
\o D^{-1}\a(p_{(1)})\c n)\\
&=&(c _{M, {}^NP}\o id_N)\circ a^{-1}_{M, {}^NP, N}\circ (id _M\o c _{N, P})(m\o(n\o p)
\end{eqnarray*}
On the another hand,
 \begin{eqnarray*}
 &&a_{{}^MN, {}^MP, M} \circ  c _{M, N\o P}\circ a_{M, N, P}((m\o n)\o p)\\
 &=&a_{{}^MN, {}^MP, M} \circ  c _{M, N\o P}(\mu(m)\o( n\o\varsigma^{-1} (p)))\\
 &=&a_{{}^MN, {}^MP, M} ((\nu\o\varsigma)(n\o \varsigma^{-1}(p))_{(0)}\o
 B^{-1}((n\o\varsigma^{-1} (p))_{(1)})\c \mu^{-1}\mu(m)\\
 &=&\nu^{2}(n_{(0)})\o( p_{(0)}\o B^{-1}\a^{-1}(\a^{-1} (p_{(1)})n_{(1)})\c \a^{-1}(m))\\
 &=&\nu^{2}(n_{(0)})\o \varsigma(\varsigma^{-1}(p)_{(0)})
 \o (B^{-1}(\varsigma^{-1}(p)_{(1)})\c(B^{-1}\a^{-1} (n_{(1)})\c \a^{-1}(m)))\\
 &=&(id _{{}^MN }\o c _{M, P})
 ((\nu^{2}(n_{(0)})\o B^{-1}(n_{(1)})\c \mu^{-1}(m))\o p)\\
 &=&(id _{{}^MN }\o c _{M, P})\circ a_{{}^MN, M, P}\ci (c _{M, N}\o id_P)
 ((m\o n)\o p)
\end{eqnarray*}

 The proof is completed.
 \hfill $\blacksquare$
\\

{\bf Lemma 3.6.} The map $c_{M,N}$
defined by $c_{M,N}(m\o n)=\nu(n_{(0)})
\o B^{-1}(n_{(1)})\c \mu^{-1}(m)$ is bijective; with inverse
 $${c}_{M,N}^{-1}(n\o m)=B^{-1}(S(n_{(1)}))\c\mu^{-1}(m) \o \nu (n_{(0)}).$$

{\bf Proof.} First, we prove $c_{M,N}{c}_{M,N}^{-1}=id$. For all $m\in M, n\in
  N$, we have
 \begin{eqnarray*}
&&c_{M,N}{c}_{M,N}^{-1}(n \o m)\\
&=&c_{M,N}(B^{-1}S(n_{(1)})\c\mu^{-1}(m)\o\nu(n_{(0)}))\\
&=&\nu(\nu(n_{(0)})_{(0)})\o B^{-1}(\nu(n_{(0)})_{(1)})\c \mu^{-1}(B^{-1}S(n_{(1)})\c\mu^{-1}(m))\\
&=&\nu^{2}(n_{(0)(0)})\o B^{-1}((n_{(0)(1)})S\a^{-1}(n_{(1)}))\c\mu^{-1}(m)\\
&=&\nu(n_{(0)})\o B^{-1}((n_{(1)1})S(n_{(1)2}))\c\mu^{-1}(m)\\
&=&\nu(n_{(0)})\o B^{-1}(\varepsilon(n_{(1)})1_{H})\c\mu^{-1}(m)\\
&=&\nu(n_{(0)})\o \varepsilon(n_{(1)})m
 = n\o m
\end{eqnarray*}

The fact that $ c_{M,N}^{-1} c_{M,N}=id$ is similar.
 This completes the proof. \hfill $\blacksquare$
\\

 Let $H$ be a monoidal Hom-Hopf algebra and
 $G={\sl Aut}_{mHH}(H)
 \times {\sl Aut}_{mHH}(H)$.
 Define $\mathcal {MHYD}(H)$ as the
  disjoint union of all $_{H}\mathcal {MHYD}^{H}(A,B)$
 with $(A,B)\in G$. If we endow $\mathcal {MHYD}(H)$
 with tensor product shown in Proposition 3.1,
 then $\mathcal {MHYD}(H)$ becomes
 a monoidal category with unit $k$.

 Define a group homomorphism
 $\,\,\varphi: G\rightarrow Aut(\mathcal {MHYD}(H)),
 \,\,\,\,\,\,\,\,
 (A, B) \,\,\mapsto \,\,\varphi(A,B)\,\,$
 on components as follows:
\begin{eqnarray*}
\varphi_{(A,B)}: {_{H}}\mathcal {MHYD}^{H}(C,D)&\rightarrow&
{_{H}}\mathcal {MHYD}^{H}((A,B)\ast
(C,D)\ast (A,B)^{-1}),\\
\quad \quad \quad \quad \quad
 \quad \varphi_{(A,B)}(N)&=& {}^{(A,B)} N,
\end{eqnarray*}
and the functor $\varphi_{(A,B)}$ acts as identity on morphisms.

 The braiding in $\mathcal {MHYD}(H)$
 is given by the family $\{c_{M,N}\}$
 in Proposition 3.5.
 So we get the following main theorem of this article.
\\

{\bf Theorem 3.7.} $\mathcal {MHYD}(H)$
is a braided $T$-category over G.
\\

We will end this paper by giving some
examples to illustrate the main theorem.

{\bf Example 3.8.}
For $\alpha=id_{H},$ we have
 $\mathcal {MHYD}(H)=\mathcal {YD}(H)$, the main constructions
 by Panaite and Staic \cite{PS2007}.
\\

{\bf Example 3.9.}   Let  $( H_{4} =k\{ 1,\ g,\ x, \ gx\,\},\a , \Delta , \varepsilon , S )$
 be the monoidal Hom-Hopf algebra (see  Example 2.5 (2)).
 Let $\mathcal {MHYD}(H_4)$ be the disjoint union of all
  categories  $_{H_4}\mathcal {MHYD}^{H_4}(A,B)$ of
  left-right Hom-$(A,B)$-Yetter-Drinfeld modules
  with ${\sl Aut}_{mHH}(H_4) \times {\sl Aut}_{mHH}(H_4).$
  Then by Example 2.5, Proposition 3.3,
  Proposition 3.5 and Theorem 3.7,
  $\mathcal {MHYD}(H_4)$ is a new braided $T$-category
  over ${\sl Aut}_{mHH}(H_4) \times {\sl Aut}_{mHH}(H_4).$
\\

Explicitly,   it is easily know  that we have a group isomorphism:
 ${\sl Aut}_{mHH}(H_{4})=\{\left(
  \begin{array}{cccc}
    1 & 0 & 0 & 0 \\
    0 & 1 & 0 & 0 \\
    0 & 0 & \lambda & 0 \\
    0 & 0 & 0 & \lambda \\
  \end{array}
\right)| 0\neq\lambda\in k\}\cong \left(k\backslash \{0\},\times\right).$
Furthermore, one has:
${\sl Aut}_{mHH}(H_{4})\times {\sl Aut}_{mHH}(H_{4})
\cong \left(k\backslash \{0\}\oplus k \backslash \{0\},\times\right).$
\\

 Let $H_{4A}\in\,\, _{H_4}\mathcal{MHYD}^{H_4}(A, id)$
and $H_{4B}\in\,\, _{H_4}\mathcal{MHYD}^{H_4}(id, B)$,
 for all $A,B\in {\sl Aut}_{mHH}(H_4).$
 Then $ (H_{4A}\otimes H_{4B},\alpha\otimes \alpha)
 \in\,\, _{H_4}\mathcal{MHYD}^{H_4}(A, B)$
 with structures as follows:
 \begin{eqnarray*}
&&h\c (x\o y)= h_{1}\c x \o h_{2}\c y,
\end{eqnarray*}
  for all $h\in H_4,\ x\in H_{4A},\ y\in H_{4B}.$
 \\

If $A=\left(
  \begin{array}{cccc}
    1 & 0 & 0 & 0 \\
    0 & 1 & 0 & 0 \\
    0 & 0 & c' & 0 \\
    0 & 0 & 0 & c' \\
  \end{array}
\right)$ and $H_{4 B}\in _{H_4}\mathcal {MHYD}^{H_4}(id,B),$
then ${}^{(A,id)}H_{4 B}=H_{4 B}$ as vector spaces, with action
$\rhd$ given by
\begin{eqnarray*}
 &&1\rhd 1=1,\ \ \ \ \ \ \ \ \ \  \ \ \ \ \ \ \ \ \
 \ \ \ \  1\rhd g = g,
 \ \ \ \ \ \ \ \ \ \ \ \ \ \ \ \ \ \ \ \  \ 1\rhd x=cx,
 \ \ \  \ \ \ 1\rhd gx = cgx\\
 &&g\rhd 1 =1, \ \ \ \ \ \ \  \ \ \ \ \ \ \ \ \
 \ \ \ \ \ \ \ g\rhd g=g,\ \ \ \ \
 \ \ \ \ \  \ \ \ \ \  \ \ \ \ \ g\rhd x=-cx,\ \ \ \  g\rhd gx=-cgx\\
 && x\rhd 1=-c(1+c'')/c'gx,
 \ \ \  x\rhd g=c(-1+c'')/c'x,
 \ \ \ \  x\rhd x=0,
 \ \ \ \ \ \ \ \  x\rhd gx=0,\\
 && gx\rhd 1=c(-1+c'')/c'gx,\
  gx\rhd g=-c(1+c'')/c'x,\ \ \
   gx\rhd x=0,\ \  \ \ \ \ \   x\rhd gx=0,
\end{eqnarray*}
and coaction $\rho_{r}$ defined by
\begin{eqnarray*}
&&\rho_{r}(1)=1\otimes1, \ \ \ \ \ \ \ \ \ \ \ \ \ \ \
\ \ \ \ \ \ \  \ \ \ \ \ \ \ \ \ \ \ \rho_{r}(g)=g\otimes g,\\
&&\rho_{r}(x)=c^{-1}(x\otimes 1)+ c^{-1}c'(g\otimes x),
 \ \ \ \ \ \ \rho_{r}(gx)=c^{-1}(gx\otimes g)+c^{-1}c'(1\otimes gx),
 \end{eqnarray*}
for all $c,c',c''\in k\backslash \{0\}$,
and ${}^{(A,id)}H_{4 B}\in  _{H_4}\mathcal {MHYD}^{H_4}(id,B).$
\\

Let $H_{4 A}\in _{H_4}\mathcal {MHYD}^{H_4}(A,id)$
and $H_{4 B}\in _{H_4}\mathcal {MHYD}^{H_4}(id,B).$
Then the braiding
$$
c_{H_{4 A},H_{4 B}}:\
H_{4 A}\otimes H_{4 B}\rightarrow {}^{(A,id)}H_{4 B}
\otimes H_{4 A}
$$
is given by
$$
c_{H_{4 A},H_{4 B}}(m\otimes n)
=\a(n_{1})\otimes n_{2}\cdot \mu^{-1}(m),
$$
for all $m\in H_{4 A}, n\in H_{4 B}.$
\\

If we consider a system of the basis of $H_{4 A}\otimes H_{4 B}$:
 $\{1\otimes 1, 1\otimes g ,
  1\otimes x,  1\otimes gx ,
  g\otimes 1 ,  g\otimes g ,
  g\otimes x ,  g\otimes gx ,
  x\otimes 1,  x\otimes g ,
  x\otimes x,  x\otimes gx ,
  gx\otimes 1 ,  gx\otimes g ,
  gx\otimes x,  gx\otimes gx \}
  $
  which is denoted by $\{e_1,e_2, \cdots , e_{16}\}$, then the braiding $c_{H_{4 A},H_{4 B}}$
  can be represented by the following matrix:
\\

$D=\left(
  \begin{array}{cccccccccccccccc}
   1 & 0 & 0 & 0 & 0 & 0 & 0 & 0 & 0 & 0 & 0 & 0 & 0 & 0 & 0 & 0 \\
    0 & 0 & 0 & 0 & 1 & 0 & 0 & 0 & 0 & 0 & 0 & 0 & 0 & 0 & 0 & 0 \\
    0 & 0 & 0 & 0 & 0 & 0 & 0 & m & 1 & 0 & 0 & 0 & 0 & 0 & 0 & 0 \\
    0 & 0 & 0 & n& 0 & 0 & 0 & 0 & 0 & 0 & 0 & 0 & 1 & 0 & 0 & 0 \\
    0& 1 & 0 & 0 & 0 & 0 & 0 & 0 & 0 & 0 & 0 & 0 & 0 & 0 & 0 & 0 \\
    0 & 0 & 0 & 0 & 0 & 1 & 0 & 0 & 0 & 0 & 0 & 0 & 0 & 0 & 0 & 0 \\
    0 & 0 & 0 & 0 & 0 & 0 & n & 0 & 0 & 1 & 0 & 0 & 0 & 0 & 0 & 0 \\
    0 & 0 & m & 0 & 0 & 0 & 0 & 0 & 0 & 0 & 0 & 0 & 0 & 1 & 0 & 0 \\
    0 & 0 & 1 & 0 & 0 & 0 & 0 & 0 & 0 & 0 & 0 & 0 & 0 & 0 & 0 & 0 \\
    0 & 0 & 0 & 0 & 0 & 0 & -1 & 0 & 0 & 0 & 0 & 0 & 0 & 0 & 0 & 0 \\
    0 & 0 & 0 & 0 & 0 & 0 & 0 & 0& 0 & 0 & 1 & 0 & 0 & 0 & 0 & 0 \\
    0 & 0 & 0 & 0 & 0 & 0 & 0 & 0 & 0 & 0 & 0 & 0 & 0 & 0 & -1 & 0 \\
    0 & 0 & 0 & 1 & 0 & 0 & 0 & 0 & 0 & 0 & 0 & 0 & 0 & 0 & 0 & 0 \\
    0 & 0 & 0 & 0 & 0 & 0 & 0 & -1 & 0 & 0 & 0 & 0 & 0 & 0 & 0 & 0 \\
    0 & 0 & 0 & 0 & 0 & 0 & 0 & 0 & 0 & 0 & 0 & 1 & 0 & 0 & 0 & 0 \\
    0 & 0 & 0 & 0 & 0 & 0 & 0 & 0 & 0 & 0 & 0 & 0 & 0 & 0 & 0 & -1 \\
  \end{array}
\right)$\\

i.e. we have

$$
c_{H_{4 A},H_{4 B}}((e_1,e_2,......e_{16})^{T})
=D(e_1,e_2,......e_{16})^{T},
$$
for all $m\neq -1,n\neq 1,$
$m,n\in k.$

  Then by Theorem 3.7,
  $\mathcal {MHYD}(H_4)$ is a new Braided $T$-category
  over $\left(k/\{0\}\oplus k/\{0\},\times\right).$
\\

\section*{ACKNOWLEDGEMENTS}

The authors would like to thank the referee for the valuable comments
 even for our English.
 This work was supported by the NSF of China (No. 11371088) and the NSF of Jiangsu Province (No. BK2012736).


\begin{thebibliography}{aa}


\bibitem{CD2006} Caenepeel S., De Lombaerde M. (2006).
 A categorical approach to Turaev's Hopf group-coalgebras.
 {\sl Commun. Algebra.} 34, 2631-2657.

 \bibitem{CG2011} Caenepeel S., Goyvaerts I. (2011).
 Monoidal Hom-Hopf algebras.
 {\sl Commun. Algebra.} 39, 2216-2240.

\bibitem{CWZ2013} Chen Y. Y., Wang Z. W., and Zhang L. Y. (2013).
  Integrals for monoidal Hom-Hopf algebras
  and their applications. {\sl J. Math. Phys.} 54, 073515.

  \bibitem{CZ2014} Chen Y. Y., Zhang L. Y. (2014).
  The category of Yetter-Drinfel'd
 Hom-modules and quantum Hom-Yang-Baxter equation.
  {\sl J. Math. Phys.} 55, 031702.


 \bibitem{FY1989} Freyd P. J., Yetter D. N. (1989). Braided
 compact closed categories with applications to low-dimensional topology.
  {\sl Adv. in Math.} 77(2), 156-182.

\bibitem{K2004} Kirillov, A. J. (2004). On $G$-equivariant
   modular categories. {\sl Math.} QA/0401119.

 \bibitem{LS2014} Liu L., Shen B. L. (2014). Radford's biproducets and
  Yetter-Drinfeld modules for monoidal Hom-Hopf algebras.
  {\sl J. Math. Phys.} 55, 031701.

  \bibitem{LW2010} Liu L., Wang, S. H. (2010). Constructing new
  braided $T$-categories over weak Hopf algebras.
  {\sl Appl. Categ. Atruct.} 18, 431-459.

 \bibitem{MS2010} Makhlouf A., Silvestrov S. D. (2010).
  Hom -algebras and Hom-coalgebras. {\sl  J. Algebra Appl.} 09,
  553-589.

  \bibitem{MS2009} Makhlouf A., Silvestrov S. D. (2007).
  Hom -Lie admissible Hom-coalgebras and Hom-Hopf algebras.
  {\sl in Generalized Lie  Theory in Mathematics, Physics and Beyond,} edited by Silvestrov S.,
  Paal E., Abramov V., Stolin A. (Springer-Verlag, Berlin, 2009), pp. 189-206;
   Preprints in Mathematical Sciences, Lund University, Centre for Mathematical Sciences,
   Centrum Scientiarum Mathematicarum, Vol. 25, 2007, LUTFMA-5091-2007 and in
    e-print arxiv:0709. 2413 [math. RA].

 \bibitem{PS2007} Panaite P., Staic M. D. (2007). Generalized (anti)
 Yetter-Drinfel'd modules as components of a braided T-category.
 {\sl Israel J. Math.} 158, 349-366.

\bibitem{S2007} Staic, M. D. (2007). A note on anti-Yetter-Drinfeld
 modules. {\sl Contemp. Math.} 441, 149-153.


\bibitem{S1969}  Sweedler, M. E. (1969).  Hopf Algebras, Benjamin, New York, 1969.


 \bibitem{T1994} Turaev V. G. (1994). Quantum Invariants of Knots and
 $3$-Manifolds. {\sl  de Gruyter Stud. Math.} 18, de Gruyter, Berlin.

 \bibitem{T2008} Turaev, V. G. (2008). Crossed group-categories.
 {\it Arab. J. Sci. Eng. Sect. C Theme Issues} 33(2C), 483-503.

\bibitem{VA2001} Virelizier, A. (2001). Alg$\grave{e}$bras de Hopf
 gradu$\acute{e}$es et fibr$\acute{e}$s plats
 sur les 3-vari$\acute{e}$t$\acute{e}$s.
 {\sl Ph. D. thesis, Universite Louis Pasteur, Strasbourg.}

\bibitem{VA2005} Virelizier, A. (2005). Involutory Hopf group-coalgebras
 and flat bundles over 3-manifolds.
 {\sl Fundam. Math.} 188, 241-270.

  \bibitem{Z2004} Zunino M. (2004). Yetter-Drinfeld modules for crossed structures.
 {\sl J. Pure Appl. Algebra.} 193, 313-343.
\end{thebibliography}
\end{document}